\documentclass[12pt]{amsart}
\usepackage{amsmath,amsthm,latexsym,amscd,amsbsy,amssymb}

\chardef\bslash=`\\ 

\makeatletter
\def\verbatim{\interlinepenalty\@M \@verbatim
  \leftskip\@totalleftmargin\advance\leftskip2pc
  \frenchspacing\@vobeyspaces \@xverbatim}
\makeatother
\newtheorem{thm}{Theorem}
\newtheorem{cor}[thm]{Corollary}

\newtheorem*{A}{Theorem A}
\newtheorem*{B}{Theorem B}
\newtheorem*{C}{Theorem C}

\theoremstyle{definition}

\theoremstyle{remark}
\newtheorem{rem}{Remark}

\numberwithin{equation}{section}

\font\f=msbm10

\begin{document}


\title[Compact Groups]
{Compact groups and absolute extensors}
\author{Alex Chigogidze}
\address{Department of Mathematics and Statistics,
University of Saskatche\-wan,
McLean Hall, 106 Wiggins Road, Saskatoon, SK, S7N 5E6,
Canada}
\email{chigogid@math.usask.ca}
\thanks{Author was partially supported by NSERC research grant.}

\begin{abstract}{We discuss compact Hausdorff groups from the point of view of the general theory of absolute extensors. In particular, we characterize the class of simple, connected and simply connected compact Lie groups as $AE(2)$-groups the third homotopy group of which is $\text{\f Z}$. This is the converse of the corresponding result of R.~Bott.}
\end{abstract}

\maketitle
\markboth{A.~Chigogidze}{Compact Groups}


The following result \cite[Theorem 6.57 and Corollary 6.59]{hofmann} plays an important role in the theory of compact groups.
\begin{A}[Structure Theorem for Compact Groups]\label{T:A} Let $G$ be a connected compact group. Then there exist a continuous homomorphism 
\[ p \colon Z_{0}(G) \times \prod\{ L_{t} \colon t \in T\} \to G\]
where $Z_{0}(G)$ stands for the identity component of the center of $G$ and $L_{t}$ is a simple, connected and simply connected compact Lie group, $t \in T$, such that $\ker (p)$ is isomorphic to a zero-dimensional central subgroup of $G$.

If $G$ itself is a compact connected $n$-dimensional Lie group, then $\ker (p)$ is finite, the indexing set $T$ is also finite and $Z_{0}(G)$ is a torus of codimension $|T|$ (i.e. $Z_{0}(G)$ is isomorphic to the product of $n-|T|$ copies of the circle group $\text{\f T}$).
\end{A} 

In many instances the above statement allows us to split a reasonably stated problem about compact connected groups into the abelian and nonabelian parts (explicitly indicated in the domain of $p$). Obviously its efficiency depends on our ability to recognize these parts separately. This is why the second part of Theorem A is much more precise than the first. One aspect of this general point of view is reflected in the following two questions:\\

\begin{itemize}
\item[(I)]
What is a (topological) characterization of tori, i.e. of (possibly uncountable) powers $\text{\f T}^{\tau}$ of the circle group $\text{\f T}$.
\item[(II)]
What is a (topological) characterization of simple, connected and simply connected compact groups?
\end{itemize}

Question (I) has recently been answered \cite[Theorem E]{chigroups} within the general theory of absolute extensors (consult \cite{book} for a comprehensive discussion of related topics):
\begin{B}\label{T:B}
The following conditions are equivalent for a compact abelian group $G$:
\begin{itemize}
\item
$G$ is a torus.
\item
$G$ is an $AE(1)$-compactum (= compact absolute extensor in dimension $1$).
\end{itemize}
\end{B} 

In these notes we give a complete solution of question (II). Characterization, as in the abelian case \cite{chigroups}, is given within the theory of $AE(n)$-spaces.

Lie groups, topologically being manifolds, are locally $n$-con\-nec\-ted for each $n \geq 0$ (even locally contractible). Two basic examples of global connectivity properties are, of course, connectedness and simple connectedness. Let us examine these concepts more closely. Connected manifolds (because of their local niceness) are arcwise connected (i.e. $0$-connected). Note that this is not true for non locally connected groups. Unlike connectedness, arcwise connectedness is a concept of a homotopy theoretical nature and fits into a general extension theory. Indeed, consider the following extension problem $E(Z,Z_{0},g_{0})$

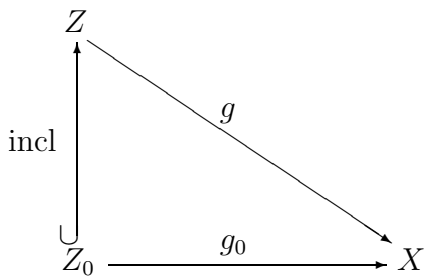
\begin{picture}(300,120)
\put(101,100){$Z$}
\put(100,10){$Z_{0}$}
\put(106,23){\vector(0,1){73}}
\put(99,20){$\cup$}
\put(80,55){$\operatorname{incl}$}
\put(118,12){\vector(1,0){105}}
\put(227,10){$X$}
\put(110,97){\vector(3,-2){115}}
\put(160,18){$g_{0}$}
\put(160,68){$g$}
\end{picture}

\noindent which asks whether a map $g_{0}$, defined on a closed subspace $Z_{0}$ of a compactum $Z$, with values in $X$, has an extension $g$ over the whole $Z$. Clearly a compactum $X$ is arcwise connected precisely when the extension problem $E(B^{1}, \partial B^{1},g_{0})$, where $B^{1}$ is the $1$-dimensional disk, is solvable for any $g_{0}$. Similarly the class of compacta for which the extension problem $E(B^{2},\partial B^{2},g_{0})$, where $B^{2}$ is the $2$-dimensional disk, is solvable for any $g_{0}$ coincides with the class of simply connected spaces. An important observation (see for instance \cite[Chapter 7, \S 53, Section IV, Theorems 1 and 1$^{\prime}$]{kura}) here is that for a {\em metrizable} locally $0$-connected compactum $X$ the $0$-connectedness of $X$ guarantees the solvability of the extension problem $E(Z,Z_{0},g_{0})$ {\em for any choice} of at most $1$-dimensional compactum $Z$, closed subspace $Z_{0}$ of $Z$, and map $g_{0}$. In other words, $X$ in this case is an absolute extensor in dimension $1$ (shortly, $AE(1)$-compactum). This fact sometimes is expressed by writing $LC^{0} \cap C^{0} = AE(1)$. These concepts can obviously be defined for any $n$ (recall that a compactum is an absolute extensor in dimension $n$, shortly $AE(n)$-compactum, if all extension problems $E(Z,Z_{0},g_{0})$ with values in $X$ and with at most $n$-dimensional $Z$ are solvable) and as is well known the strictly decreasing sequences
\[ AE(1) \supset AE(2) \supset \cdots \supset AE(n) \supset AE(n+1) \supset \cdots \]
and 
\[ LC^{0}\cap C^{0} \supset LC^{1}\cap C^{1} \supset \cdots \supset LC^{n-1}\cap C^{n-1} \supset LC^{n}\cap C^{n} \supset \cdots \]
are identical (i.e., $AE(n) = LC^{n-1}\cap C^{n-1}$) in the presence of metrizability. Nevertheless these concepts differ (i.e. $AE(n) \subsetneq LC^{n-1}\cap C^{n-1}$) in general. For compact groups the difference between these concepts is very delicate. To see this it suffices to compare Theorem B with the result from \cite{shelah}: the validity of the statement ``an arcwise connected (i.e. $0$-connected) compact abelian group is a torus" is undecidable within ZFC.

Thus every connected Lie group topologically is an  $AE(1)$-space and every connected and simply connected Lie group is an $AE(2)$-space.  

\begin{C}\label{T:C}
The following conditions are equivalent for a compact group $G$:
\begin{itemize}
\item[(a)]
$G$ is a simply connected $AE(1)$-compactum.
\item[(b)]
$G$ is an $AE(2)$-compactum.
\item[(c)]
$G$ is an $AE(3)$-compactum.
\item[(d)]
$G$ is a product of simple, connected and simply connected compact Lie groups.
\end{itemize}
\end{C}
\begin{proof}
Implications (c) $\Longrightarrow$ (b) and (b) $\Longrightarrow$ (a) are trivial.

(a) $\Longrightarrow$ (d).
It follows from the proof of Theorem B \cite{hofmann} that there exists a continuous homomorphism $\alpha \colon G^{\prime} \times Z_{0}(G) \to G$ where $G^{\prime}$ is the commutator subgroup of $G$, $Z_{0}(G)$ is the connected component of the center of $G$, and $\ker (\alpha )$ is a zero-dimensional compact group isomorphic to the (closed and central) subgroup $G^{\prime}\cap Z_{0}(G)$ of $G$. Further, there exists a continuous homomorphism $\beta \colon \prod\{ L_{t} \colon t \in T\} \to G^{\prime}$ where each $L_{t}$, $t \in T$, is a simple, connected and simply connected compact Lie group, and $\ker (\beta )$ is zero-dimensional. Obviously, the kernel of the homomorphism $p$, defined as the composition 
\[ p = \alpha \circ (\beta \times \operatorname{id}_{Z_{0}(G)}) \colon \prod\{ L_{t} \colon t \in T\} \times Z_{0}(G) \xrightarrow{\beta \times \operatorname{id}} G^{\prime}\times Z_{0}(G) \xrightarrow{\alpha} G ,\]

\noindent is zero-dimensional. Let us show that $|\ker (p)| = 1$. Indeed, assuming that $|\ker (p)| > 1$ and remembering that $\dim\ker (p) = 0$, we can find a compact group $\tilde{G}$ and two continuous homomorphisms 
\[ p_{1} \colon \prod\{ L_{t} \colon t \in T\} \times Z_{0}(G) \to \tilde{G}\;\;\text{and}\;\; p_{2} \colon \tilde{G} \to G\]
such that $p = p_{2}p_{1}$, $\ker (p_{2})$ is finite and $|\ker (p_{2})| >1$. Clearly $p_{2}$ is a bundle (since $\ker (p_{2})$ is a finite group) and $G$, as an $AE(1)$-compactum, is arcwise connected and locally arcwise connected. Also, by (a), $\pi_{1}(G) = 0$. Consequently, by \cite[Corollary on p.66]{steenrod}, $p_{2}$ is a trivial bundle, i.e. $\tilde{G} \cong G \times \ker (p_{2})$. Now observe that $\tilde{G}$, as a continuous image of $\prod\{ L_{t} \colon t \in T\} \times Z_{0}(G)$, is connected. This implies that $\ker (p_{2})$ is a singleton which is impossible. This contradiction shows that $|\ker (p)| = 1$ and hence $\prod\{ L_{t} \colon t \in T\} \times Z_{0}(G) \cong G$. From this we conclude that $Z_{0}(G)$, as a retract of $G$, is an $AE(1)$-compactum with $\pi_{1}(Z_{0}(G)) = 0$. An abelian group which is an $AE(1)$-compactum is a torus (Theorem B). But the only simply connected torus is the singleton. Thus $|Z_{0}(G)| = 1$ and $G \cong \prod\{ L_{t} \colon t \in T\}$. 

(d) $\Longrightarrow$ (c). Suppose that $G$ is a product $\prod\{ L_{t} \colon t \in T\}$ of simple, connected and simply connected compact Lie groups. By Theorem of H.~Cartan \cite[Theorem 3.7]{mimura}, $\pi_{2}(L_{t}) = 0$, $t \in T$. By the local contractibility of $L_{t}$, this means that $L_{t}$ is an $AE(3)$-compactum, $t \in T$. Then $G$, as a product of $AE(3)$-compacta, is an $AE(3)$-compactum. The proof is complete.
\end{proof}
\begin{rem}
Implication (a) $\Longrightarrow$ (b) of Theorem C extends a result of H.~Cartan which states that if $G$ is a connected and (semi-) simple compact Lie group, then $\pi_{2}(G) = 0$.
\end{rem}
We now present corollaries of Theorem C. The following answers question (II).

\begin{cor}\label{C:bott}
The following conditions are equivalent for a non-trivial compact group $G$:
\begin{itemize}
\item[(i)]
$G$ is an $AE(2)$-group with $\pi_{3}(G) = \text{\f Z}$.
\item[(ii)]
$G$ is a simple, connected and simply connected Lie group.
\end{itemize}
\end{cor}
\begin{proof}
(i) $\Longrightarrow$ (ii). By Theorem C, $G$ is isomorphic to the product $\prod\{ L_{t} \colon t \in T\}$ of simple, connected and simply connected compact Lie groups. Let $T^{\prime} = \{ t \in T \colon |L_{t}| > 1\}$. By (i) and a result of R.~Bott \cite[Theorem 3.8]{mimura}, we then have
\[ \text{\f Z} = \pi_{3}(G) = \pi_{3}\left(\prod\{ L_{t} \colon t \in T^{\prime} \}\right) = \prod\{ \pi_{3}(L_{t}) \colon t \in T^{\prime}\} = \text{\f Z}^{T^{\prime}} .\]
This obviously is possible precisely when $|T^{\prime}| = 1$ which implies that $G$ is isomorphic to some $L_{t}$ and hence is a simple, connected and simply connected compact Lie group.

(ii) $\Longrightarrow$ (i). Every connected and simply connected Lie group is an $AE(2)$-space. Since $G$ is simply connected, it is non abelian. Since $G$ is simple, by the above cited result of R.~Bott, $\pi_{3}(G) = \text{\f Z}$.
\end{proof}

\begin{rem}
Implication (ii) $\Longrightarrow$ (i) of Corollary \ref{C:bott} is due to R.~Bott.
\end{rem}

\begin{cor}\label{C:nae4}
There is no non-trivial compact $AE(4)$-group.
\end{cor}
\begin{proof}
Let $G$ be a compact $AE(4)$-group. Since every $AE(4)$-compactum is an $AE(3)$-compactum, it follows from Theorem C, that $G$ is a product $\prod\{ L_{t} \colon t \in T\}$ of simple, connected and simply connected compact Lie groups. Since each retract of an $AE(4)$-compactum is an $AE(4)$-compactum, it follows that $L_{t}$ is an $AE(4)$-compactum, $t \in T$. Then $\pi_{3}(L_{t}) = 0$ for each $t \in T$. This contradicts Corollary \ref{C:bott} and completes the proof.
\end{proof}

The concept of an $ANE(4)$-space (=absolute {\em neighborhood} extensor in dimension $4$) is obtained by requesting (in the extension problem $E(Z,Z_{0},g_{0})$ with any at most $4$-dimensional space $Z$, any closed subspace $Z_{0}$ of $Z$, and any continuous $g_{0}$) the existence of an extension $g$ (of $g_{0}$) defined not on the whole $Z$ (as we did while defining the concept of $AE(4)$) but only on a neighborhood of $Z_{0}$ in $Z$. Clearly every $ANE(4)$-space is locally $3$-connected but not vice versa. It is well known that a locally compact group is a Lie group if and only if it is finite-dimensional and locally connected \cite{pontr}. It is easy to produce examples of non-metrizable locally connected compact groups, whereas a finite dimensional connected locally compact group is always metrizable \cite{skl}. Roughly speaking our next Corollary shows what degree of local connectedness forces an infinite dimensional compact group to be metrizable

\begin{cor}\label{C:ane4}
Every locally compact $ANE(4)$-group is metrizable.
\end{cor}
\begin{proof}
Let $G$ be a locally compact $ANE(4)$-group. By \cite[Theorem 16$^{\prime}$]{skl}, $G$ is homeomorphic to the product $\text{\f Z} \times \left(\text{\f Z}_{2}\right)^{\tau} \times \text{\f R}^{n} \times L$ for some integer $n$, a cardinal number $\tau$ and a compact group $L$. Since each retract of an $ANE(4)$-space is an $ANE(4)$-space, we conclude that $\tau$ is finite. Therefore it suffices to show that $L$ is metrizable.

By \cite[Lemma 8.2.1]{book}, $L$ is isomorphic to the limit of an $\omega$-spectrum ${\mathcal S}_{L} = \{ L_{\alpha}, p_{\alpha}^{\beta}, A\}$ consisting of metrizable compact groups $L_{\alpha}$ and continuous homomorphisms $p_{\alpha}^{\beta} \colon L_{\beta} \to L_{\alpha}$. Observe that $L$, as a retract of $G$, is an $ANE(4)$-compactum. Consequently, by \cite[Theorems 1.3.4 and 6.3.2]{book}, there exists a cofinal (in particular, non-empty) subset $B$ of the indexing set $A$ such that the limit projection $p_{\alpha} \colon L \to L_{\alpha}$ is $4$-soft for each $\alpha \in B$. Then, by \cite[Corollary 6.1.17 and Proposition 6.1.21]{book}, the fiber $\ker (p_{\alpha})$ is a compact $AE(4)$-group. By Corollary \ref{C:nae4}, $|\ker (p_{\alpha})| =1$. Therefore $p_{\alpha}$ is an isomorphism. Since $L_{\alpha}$ is metrizable, so is $L$.
\end{proof}


\end{document}